
\input amstex
\documentstyle{amsppt}
\magnification=1200
\hsize=16.5truecm
\vsize=23.3truecm

\catcode`\@=11
\redefine\logo@{}
\catcode`\@=13

\define\vth{\vartheta}

\define\ml{\operatorname{Exp-Lift}}

\define \bz{\Bbb Z}
\define \bq{\Bbb Q}
\define \br{\Bbb R}
\define \bc{\Bbb C}
\define \bh{\Bbb H}

\define\pd#1#2{\dfrac{\partial#1}{\partial#2}}

\define \Mn{M^{[n]}}

\TagsOnRight
\NoBlackBoxes
\document

\topmatter
\title
Complex vector bundles and Jacobi forms
\endtitle

\rightheadtext
{Complex vector bundles and Jacobi forms}
\leftheadtext{V. Gritsenko}

\author
V.  Gritsenko 
\endauthor

\thanks{To appear in the Proceedings  of the  Symposium 
``Automorphic forms and L-functions", RIMS, Kyoto, Japan,
January 25--29,  1999}
\endthanks

\address
D\'epartement de Math\'ematique Universit\'e Lille I
\newline
${}\hskip 9pt $ 59655 Villeneuve d'Asq Cedex, France
\endaddress

\abstract
The elliptic genus (EG) of a compact complex  manifold was introduced
as a holomorphic Euler characteristic of some formal power series
with vector bundle coefficients.  EG is an automorphic form in two variables 
only if the manifold  is a Calabi--Yau manifold.
In physics such a function appears as the partition function of
$N=2$ superconformal field theories.
In these notes  we define  the modified Witten genus or  
the automorphic correction of elliptic genus. 
It  is  an automorphic function  in two variables  for  an arbitrary 
holomorphic vector bundle  over a compact  complex manifold.
This paper is an exposition of the talks given by the author
at Symposium ``Automorphic forms and L-functions" at RIMS, Kyoto
(January, 27, 1999) and at Arbeitstagung in Bonn (June, 20, 1999).
\endabstract

\endtopmatter

\document
\document

\head
Introduction 
\endhead

In these notes we present a  link between the theory of automorphic forms
and geometry.
For an arbitrary compact spin manifold one can define its elliptic genus.
It is a modular form in one variable with respect to a congruence
subgroup of level $2$ (see, for example, \cite{W1},
\cite{L}, \cite{HBJ}).
For  a   compact complex manifold one can define its elliptic genus
$\phi(M;\tau,z)$
as a function in two complex variables 
(see, for example, \cite{W2}, \cite{EOTY},  \cite{H\"o},
\cite{FY},  \cite{KYY}).
In the last case the elliptic genus is 
the holomorphic Euler characteristic of a formal power series 
with vector bundle coefficients.
If the first Chern class $c_1(M)$ of the complex manifold is equal to zero
in $H^2(M,\br)$, then the elliptic genus is a  weak Jacobi modular form
(with integral Fourier coefficients) of weight $0$ and index $d/2$,
where $d=$dim${}_\bc(M)$. 
The same modular form appears in physic as the partition function 
of $N=2$ super-symmetric sigma model whose target space is 
the given Calabi--Yau manifold.
We note that any ``good" partition function has appeared in physic 
is an automorphic form with respect to some   group. 
This reflects the fact  that   physical models have some additional 
symmetries.
If $c_1(M)\ne 0$, then the elliptic genus
$\phi(M;\tau,z)$  of $M$ is not an automorphic form.
In these notes  we define {\it the  modified Witten genus}  (MWG) or  
{\it the automorphic correction of elliptic genus}
of  an arbitrary holomorphic vector bundle over a compact complex manifold
and we  briefly study its properties. 
This new object is always an automorphic form in two variables.

We are going to present here  automorphic aspects of the theory. 
In the proof of  the theorem that the modified 
Witten genus is a Jacobi form we use a nice formula
which relates the Jacobi theta-series, its logarithmic derivative,
the quasi-modular Eisenstein series $G_2(\tau)$ and all the  derivatives
of the Weierstrass $\wp$-function (see Proposition 1.4 bellow).
The fact that MWG of a  vector bundle is a Jacobi modular form
has a number of applications in the theory of vector bundles.
Among them we have relations between $\hat A$-genera and  other
cohomological invariants. One can get from this  some obstructions  if 
the rank of the vector bundle is smaller that the dimension
of the manifold. The special values of the MWG have also  some nice 
arithmetic properties.

One can see from the definition that the coefficients of MWG
are rationals with bounded denominators.
To get applications  to the theory
of complex manifolds  we study the  $\bz$-structure of the bigraded ring 
$$
J^{w, \bz}_{*,*/2}=\bigoplus_{k\in \bz, \,m\in 
\frac {1}2\bz}J^{w,\bz}_{k,m}
$$
of all weak  Jacobi forms of all weights and indices 
with integral Fourier  coefficients. 
We prove that this  ring has $14$ generators over $\bz$ (see Theorem 2.3):
$$
J_{*, */2}^{w, \bz}=\bz[E_4, E_6, \Delta,
E_{4,1}, E_{4,2}, E_{4,3}, 
E_{6,1}, E_{6,2}, E'_{6,3}, 
\phi_{0,1}, \phi_{0,\frac {3}{2}}, \phi_{0,2},  \phi_{0,4},
\phi_{-1, \frac{1}2}],
$$
where $E_4$, $E_6$ and $\Delta$ are $SL_2(\bz)$-generators,
$E_{4,1},\dots E_{6,3}$ are the Eisenstein--Jacobi series
and  $E_{6,3}'$ is a modification of $E_{6,3}$.
The four functions $\phi_{0,1}, \dots, \phi_{0,4}$  of weight $0$
generate the graded $\bz$-ring  $J^{w, \bz}_{0,*/2}$ of 
the   Jacobi forms of weight $0$ with integral Fourier coefficients.
These Jacobi forms  are related to Calabi--Yau manifolds of dimension 
$d=2,\,3,\,4,\,8$. The Jacobi forms  $\phi_{0,1}$,
$\phi_{0,2}$, $\phi_{0,3}=\phi_{0,3/2}^2$  
are the  generating functions for the multiplicities
of  positive roots of the three main  generalized Lorentzian Kac--Moody
Lie algebras of Borcherds type constructed in \cite{GN1--GN4}
(see  also \S 4 of this paper). The form $\phi_{0,4}$ determines
the "most odd" even Siegel theta-constant.
(see \cite{G4, \thetag{2.11}}).

The $q^0$-term of the Fourier expansion ($q=e^{2\pi i \tau}$) 
of the elliptic genus  $\phi(M;\tau,z)$
is essentially equal to the Hirzebruch 
$\chi_y$-genus of the manifold.
Thus we can analyze the arithmetic properties of the  
$\chi_y$-genus of the complex manifold with $c_1(M)=0$ and 
its special values  such as the signature ($y=1$) and the Euler number 
($y=-1$) in  terms of    Jacobi  forms.
For example, we prove that  the   Euler number
of a Calabi-Yau manifold  $M_d$ of dimension $d$ satisfies 
$$
e(M_d)\equiv 0\mod 8\qquad \text{if}
\quad d\equiv 2 \mod 8
$$
(see Proposition 2.6).
The special values of the generators of the Jacobi ring 
at $z=\frac {1}2$, $\frac {1}3$,
$\frac {1}4$ are related to Hauptmoduls of the fields
of modular functions. Using this fact we prove that 
$$
\chi_{y=\zeta_3}(M_{d})\equiv 0\mod 9\qquad \text{if}
\quad d\equiv 2 \mod 6
$$
(see Proposition 2.7).
The special values of MWG are related with $\hat A$-genera.
Some other constructions (for example, $\hat A_2^{(2)}$-genus, 
the second quantized elliptic genus) and  other applications 
to the theory of vector bundles one can find in my course of lectures
given at RIMS, Kyoto University, at our joint seminar with K. Saito
in 1998. I would like to take this opportunity 
to express  my   gratitude to all members of K. Saito's seminar.

\head
\S 1. Automorphic correction of elliptic genus
\endhead

Let $M$ be a complex compact manifold $M$ of (complex) dimension
$d$ and let $E$ be 
a complex  vector bundle over $M$.
Let us fix  two  formal variables
$q=\exp(2\pi i \tau)$ and  $y=\exp(2\pi i z)$,
where $\tau\in \bh_1$ (the upper half-plane) and 
$z\in \bc$.  One defines a formal power series
${\bold E}_{q,y}\in K(M)[[q, y^{\pm1}]]$ 
$$
{\bold E}_{q,y}=  \bigotimes_{n= 0}^{\infty}
{\bigwedge}_{-y^{-1}q^{n}}E^*
\otimes 
 \bigotimes_{n= 1}^{\infty}{\bigwedge}_{-y q^n}  E\otimes 
 \bigotimes_{n= 1}^{\infty} S_{q^n} T_M^*
\otimes 
 \bigotimes_{n= 1}^{\infty} S_{q^n} T_M ,
\tag{1.1}
$$
where $T_M$ denotes the holomorphic tangent bundle of $M$ and 
$$ 
{\bigwedge}_x E=\sum_{k\ge 0}\,  (\wedge^k E) x^k , \qquad
S_x E= \sum_{k\ge 0} (S^k E) x^k
$$
are formal power series with   exterior powers and symmetric powers 
of a bundle $E$ as coefficients.
We propose the following 

\definition{Definition 1.1}
{\it Modified Witten genus} (MWG)  
of a complex  vector bundle $E$ of rank $r$
over a compact complex manifold $M$ of dimension $d$
is defined as follows
$$
\align
{}&\chi(M, E;\tau,z)=
q^{(r-d)/12}y^{r/2}\int_M  
\exp \biggl(\frac 1{2}\bigl(c_1(E)-c_1(T_M)\bigr)\biggr)\cdot\\
\vspace{2\jot}
{}&
\exp\biggl(\bigl(p_1(E)-p_1(T_M)\bigr)\cdot G_2(\tau)\biggr)
\exp\biggl(-\frac  {c_1(E)}{2\pi i} 
\frac {\vth_z(\tau,z)}{\vth(\tau,z)}\biggr)
\,\hbox{ch}({\bold E}_{q,y})\, \hbox {td}(T_M),
\endalign
$$
where 
$c_1(E)$ and $p_1(E)$ are the first Chern and Pontryagin class
of $E$, 
$\hbox {td\,}$ is the Todd class, $\hbox{ch}({\bold E}_{q,y})$
is the Chern  character applied to each  coefficient
of the formal power series  and 
the integral $\int_M$ denotes the evaluation of 
the top degree differential form
on the fundamental cycle of the manifold.
\enddefinition

In the definition we use Jacobi theta-series of level two 
$\vartheta(\tau,z)=-i\vartheta_{11}(\tau, z)$:
$$
\vartheta(\tau ,z)=\hskip-2pt\sum\Sb n\equiv 1\, mod\, 2 \endSb
\,(-1)^{\frac{n-1}2}
q^{\frac {n^2}{8}} y^{\frac {n}{2}}=
-q^{1/8}y^{-1/2}\prod_{n\ge 1}\,(1-q^{n-1} y)(1-q^n y^{-1})(1-q^n),
$$
$\vartheta_z(\tau,z)=\pd {\vth}{z}(\tau,z)$ and 
$
G_2(\tau)=-\frac 1{24}+\sum_{n=1}^{\infty}\sigma_1(n)\,q^n
$
is a quasi-modular  Eisenstein series of weight $2$,
where $\sigma_1(n)$ is the sum of all positive divisors of $n$.

\example{1.2. Witten genus}
As a particular   case of the definition above one obtains
the Witten genus (see \cite{W1}, \cite{W2}, \cite{L}, \cite{HBJ}).
Let assume that $M$ admits a spin structure
(i.e., the second Witney-Stiefel class $w_2(M)$ is zero  
or $c_1(T_M)\equiv 0 \mod 2$)  and $p_1(M)=0$.
Let $E=\bc^r$ be the trivial vector bundle of rank $r$ over $M$. 
Then
$\hbox{ch\,}(\bigwedge_x E)=(1+x)^r$ and
$$
q^{r/12}y^{r/2}
\hbox{ch\,}\biggl(\bigotimes_{n= 0}^{\infty}
{\bigwedge}_{-y^{-1}q^{n}}E^*
\otimes 
 \bigotimes_{n= 1}^{\infty}{\bigwedge}_{-y q^n}  E
\biggr)=\biggl(\frac{\vth(\tau,z)}{\eta(\tau)}\biggr)^r.  
$$
Thus
$$
\multline
q^{d/{12}}\chi(M, \bc^r;\tau,z)=
\frac{\vth(\tau,z)^r}{\eta(\tau)^r}\int_M
\prod_{i=1}^d 
\frac {x_i/2}{\hbox{sinh}(x_i/2)}
\prod_{n=1}^{\infty}\frac {1}{(1-q^ne^{x_i})(1-q^ne^{-x_i})}\\
=
\hat A\bigl(M, \bigotimes_{n=1}^{\infty} S_{q^n} 
(T_M\oplus {T_M}^*)\bigr)\frac{\vth(\tau,z)^r}{\eta(\tau)^{r}}=
\hbox{Witten genus\,}(M) \frac{\vth(\tau,z)^r}{\eta(\tau)^{r+2d}}.
\endmultline
$$
If we take the trivial vector bundle of rank $0$, then
$$
\chi(M, 0;\tau,z)=\chi(M;\tau)=
\frac{\hbox{Witten genus\,}(M)}{\eta(\tau)^{{2d}}}.
$$
This is an automorphic function  in $\tau$ with respect to $SL_2(\bz)$.
\endexample

\example{1.3. Elliptic genus of a Calabi--Yau manifold}
This case is of some interest in physics.
Let $E=T_M$ and $c_1(T_M)=0$. Then there are no correction terms
of type $\exp(\dots)$ in Definition 1.1. 
Thus the MEG of $T_M$ is, up to the factor  $y^{d/2}$, 
the Euler--Poincar\'e characteristic of the 
element ${\bold E}_{q,y}$. This function is called
{\it elliptic genus} of the Calabi--Yau manifold $M$
or genus one partition function of the super-symmetric 
$(2,2)$-sigma model whose target space is $M$:
$$
\chi(M, T_M;\tau,z)=
\text{Elliptic genus}\,(M;\tau,z)=
y^{d/2}\int_M  
\,\hbox{ch}({\bold E}_{q,y})\, \hbox {td}(T_M). 
$$
According to  the Riemann--Roch--Hirzebruch theorem one can see
that the $q^0$-term of 
$\chi(M;\, \tau,z)$   is essentially
the Hirzebruch $\chi_y$-genus of the manifold $M$
(with interchanging of $y$ with $-y$):
$$
\gather
\chi(M;\, \tau,z)=\sum_{p=0}^d (-1)^p\chi_p(M)\, y^{\frac{d}2-p}+
\tag{1.2}
\\
q\biggl(\sum_{p=-1}^{d+1}
(-1)^{p}y^{-p}
\bigl(\chi_p(M, T_M^*)-\chi_{p-1}(M, T_M^*)+
\chi_p(M, T_M)-\chi_{p+1}(M, T_M)\bigr)\biggr)
+\dots
\endgather
$$
where 
$
\chi(M, E)=\sum_{q=0}^d (-1)^q \,\hbox{dim}\, H^q (M, E)
$ 
and 
$
\chi_p(M, E)=\chi(M,\, \wedge^p T_M^*\otimes E)
$
or, for a K\"ahler manifold,
$\chi_p(M)=\sum_{q}(-1)^{q}h^{p,q}(M)$.
We remark that in this case  every  Fourier coefficient of 
the elliptic genus is equal to the index of the  Dirac operator
twisted with a corresponding  vector bundle coefficient of the formal
power series ${\bold E}_{q,y}$.

It is known that  the elliptic genus of a Calabi--Yau manifold 
is a modular form 
in variables $\tau$ and $z$ (see \cite{H\"o}, \cite{KYY}, \cite{Li}), 
i.e.,  it is  a weak Jacobi form  of weight $0$ and index $d/2$.
{\it If $c_1(T_M)\ne 0$, then the elliptic genus
of $M$ defined above  is not a modular form}.
We add the three correction  factors in Definition 1.1  in order
to obtain a function with a good behavior with respect to the modular 
transformations in  $\tau$ and $z$.
If $E=T_M$ and $c_1(T_M)\ne 0$, then 
the integral in Definition 1.1 contains  the only 
correction term 
$$
\exp\biggl(-\frac {c_1(T_M)}{2\pi i} \frac {\vth_z}{\vth}(\tau,z)\biggr).
$$
Thus the elliptic genus of $M$ (as a function in  two variables) is equal 
to the zeroth term in a sum
of $d+1$ summands of the modified  genus. These summands  correspond 
to all powers of the first Chern class of $M$
$$
\chi(M, T_M;\tau,z)=\text{Elliptic genus}\,(M;\tau,z)
+\sum_{n=1}^d \bigl(\int_M c_1(M)^n (\dots)\bigr).  
$$
In general the elliptic genus is not an automorphic form 
but the modified elliptic genus is.
The main result of this section is
\endexample

\proclaim{Theorem 1.2} Let $E$ be a holomorphic  vector bundle 
of rank $r$ over a compact complex  manifold $M$ of dimension
$d$. Let $\chi(M,E;\tau, z)$ be  the modified Witten genus. 
Then the product
$$
\chi(M,E;\tau, z)\,\biggl(\frac{\vth(\tau,z)}{\eta(\tau)}\biggr)^{d-r}
$$
is a weak  Jacobi form of weight $0$ and index $d/2$.
\endproclaim

\remark{Remark}
We note that the definition of MWG, Theorem 1 and its proof are valid
for an arbitrary compact $Spin^{c}$-manifold $M$ and a Hermitian 
$C^{\infty}$-vector bundle $E$ over $M$. In particular, we have the same
result for an almost complex manifolds.
\endremark

\smallskip
First we recall the definition of Jacobi forms of the type we need
in this paper.
Let $t\ge 0$ and $k$   be integral or half-integral. 
Let  $v$ be a character of finite order
(or a multiplier system for half-integral  $k$) of $SL_2(\bz)$.
A holomorphic  function
$\phi(\tau ,z)$ on $\bh_1\times \bc$ is called
a  {\it  weak Jacobi form of
weight $k$ and index $t$ with  character} $v$
if it  satisfies the functional equations
$$
\phi\biggl(\frac{a\tau+b}{c\tau+d},\,\frac z{c\tau+d}\biggr)=
v(\gamma)(c\tau+d)^k \,e^{2\pi i t\tsize{\frac { c z^2}{c\tau+d}}}
\,\phi(\tau,z)
\qquad (
\gamma=\left(\smallmatrix a&b\\c&d\endpmatrix \in SL_2(\bz))
\tag{1.3a}
$$
and 
$$
\phi(\tau, z+\lambda \tau+ \mu)=
(-1)^{2t(\lambda+\mu)}
\,e^{-2\pi i t(\lambda^2 \tau+2\lambda z)}\phi(\tau,z)
\qquad (\lambda, \mu\in \bz) \tag{1.3b}
$$
and $\phi(\tau ,z)$ has  the Fourier expansion of the type
$$
\phi(\tau ,z)=\sum\Sb n\ge 0\ l \in t+\bz\\
\vspace{0.5\jot} \endSb
f(n,l)\, q^n y^l.
$$
We denote the space of all week Jacobi forms of weight $k$,
index $t$ and character (or multiplier system) $v$ by $J_{k,t}^w(v)$.
The space $J_{k,t}^w(v)$ is finite dimensional (see \cite{EZ}).
The only difference with \cite{EZ} is that
we admit Jacobi forms of half-integral index.

\proclaim{Corollary 1.3} In the conditions of Theorem 1.2 we have that
$$
\chi(M,E;\tau,z)\in J_{0,\frac{r}2}^w(v_\eta^{2(r-d)})
\qquad{r\ge d}
$$
or 
$$
\chi(M,E;\tau,z)\in J_{0,\frac{r}2}^{mer}(v_\eta^{2(r-d)})
\qquad{r< d}
$$
with a possible pole of order not bigger than $d-r$ along $z=0$.
\endproclaim
\demo{Proof of Theorem 1.2}To prove  the theorem we represent 
$\chi(M,E;\tau, z)$  in terms of the theta-series.
Let $c(E)$ be the total Chern class of the vector bundle  $E$
$$
c(E)=\sum_{i=0}^r c_i(E)=\prod_{i=1}^{r}(1+x_i)
$$
where $x_i=2\pi i \xi_i$ ($1\le i\le r$) are the formal
Chern roots  of $E$. We denote by  $x_j'=2\pi i \zeta_j$ ($1\le j\le d$) 
the  Chern roots  of $T_M$.
We recall that
$$
\hbox{ch\,}({\bigwedge}_t E)=\prod_{i=1}^r(1+te^{x_i}),
\qquad
\hbox{ch\,}(S_t E)=\prod_{i=1}^r \frac{1}{1-te^{x_i}}.
$$
According to the last formulae we have
$$
\hbox{ch\,}({\bold E}_{q,y})\,\hbox{td\,}(T_M)=
 \prod_{n=1}^{\infty}\, \prod_{j=1}^d\prod_{i=1}^r\, 
\frac{(1-q^{n-1}y^{-1}e^{-x_i})(1-q^{n}ye^{x_i})}
{(1-q^{n-1}e^{-x_j'})(1-q^{n}e^{x_j'})}x_j'\,.
$$
Therefore
$$
\multline
q^{(r-d)/12}y^{r/2} 
\exp \biggl(\frac 1{2}(c_1(E)-c_1(T_M)\biggr)
\hbox{ch\,}({\bold E}_{q,y})\hbox{td\,}(T_M)=\\
(-1)^{r-d}\prod_{i=1}^r \frac{\vth(\tau, -z-\xi_i)}{\eta(\tau)}
\prod_{j=1}^d \frac{\eta(\tau)}{\vth(\tau, -\zeta_j)} (2\pi i \zeta_j).
\endmultline
\tag{1.4}
$$
Putting  the last expression  under the integral 
we obtain the following formula for 
the modified Witten genus
$$
\align
\chi(M,E;\tau, z)=
\int_M&
\,\prod_{i=1}^r \exp\biggl(-4\pi^2 G_2(\tau)\xi_i^2
-\frac{\vth_z}{\vth}(\tau,z
)\xi_i\biggr)
\,\frac {\vth(\tau, z+\xi_i)}{\eta(\tau)}\times\\
{}&\ \prod_{j=1}^d\exp\bigl(4\pi^2 G_2(\tau)\zeta_i^2\bigr)\,
\frac{\eta(\tau)}{\vth(\tau,\zeta_j)}(2\pi i \zeta_j).
\tag{1.5}
\endalign
$$
We shall calculate the top differential form under the integral using
Proposition 1.4 bellow.
In order to formulate it we need to recall the definition of 
the Weierstrass $\wp$-function 
$$
\wp(\tau,z)=z^{-2}+\sum\Sb\omega\in \bz\tau+\bz\\ \omega\ne 0\endSb
\bigl((z+\omega)^{-2}-\omega^{-2}\bigr)
\in J_{2,0}^{mer}
\tag{1.6}
$$
which is a meromorphic Jacobi form of weight $2$ and index $0$ 
with pole of order  $2$ along $z\in \bz \tau+\bz$.
\enddemo

\proclaim{Proposition 1.4}The following formula is valid 
$$
\exp\biggl(-4\pi^2 G_2(\tau)\xi^2
-\frac{\vth_z}{\vth}(\tau,z)\xi\biggr)
{\vth(\tau, z+\xi)}=
\vth(\tau, z)\exp\biggl(-\sum_{n\ge 2} \wp^{(n-2)}(\tau,z) 
\frac{\xi^n}{n!}\biggr)
$$
where 
$
\wp^{(n)}(\tau,z)=\dfrac{\partial^n}{\partial z^n}\wp(\tau,z).
$
\endproclaim
\demo{Proof}
One of the main examples of weak Jacobi forms of half-integral index
with trivial $SL_2$-character is  
$$
\phi_{-1,1/2}(\tau, z)={\vth(\tau,z)}/{\eta(\tau)^3}
=(y^{1/2}-y^{-1/2})+q(\dots)\in J_{-1,\frac{1}2}^w. 
$$

The Jacobi form  $\phi_{-1, \frac 1{2}}$  has 
the following exponential representation as a Weierstrass $\sigma$-function
(see, for example,  review \cite{Sk})
$$
\phi_{-1,\frac{1}2}(\tau,z)=\frac {\vth(\tau,z)}{\eta(\tau)^3}=
(2\pi i z)
\exp\biggl(
-
\sum_{k\ge 1} \frac {2}{(2k)!}G_{2k}(\tau)(2\pi i z)^{2k} \biggr)
\tag{1.7}
$$
where
$
G_{2k}(\tau)=-{B_{2k}}/{4k}+\sum_{n=1}^{\infty}\sigma_{2k-1}(n)q^n
$
is the  Eisenstein series of weight $2k$.
(For each $\tau\in \bh_1$ the product is normally convergent in $z\in \bc$.)
Since one can obtain the Weierstrass $\wp$-function as the second derivative
of the Jacobi theta-series
$
\dfrac{\partial^2}{\partial z^2}\log \vth(\tau,z)=
-\wp(\tau,z)+8\pi^2 G_2(\tau),
$
the identity \thetag{1.7} implies that 
$$
\wp^{(n-2)}(\tau,z)=\frac{(-1)^n(n-1)!}{z^{n}}
+2\sum_{k\ge 2,\, 2k\ge n}
(2\pi i)^{2k} G_{2k}(\tau)\frac {z^{2k-n}}{(2k-n)!}.
$$
After that  the formula of the lemma follows by direct calculation.
\enddemo

Now we can finish the proof of Theorem 1.2.
According to the formula of Proposition 1.4   we can split
the Chern roots $x_i$ ($1\le i\le r$)
of the vector bundle $E$ and the Chern roots 
$x'_j$ ($1\le j\le d$) of  the manifold $M$
under  the integral in \thetag{1.5}, i.e.,
$$
\chi(M,E;\tau, z)=
\frac {\vth^r}{\eta^{r+2d}} \int_M  P(E;\tau,z)\cdot W(M;\tau).
\tag{1.8}
$$
The first factor depends only on the vector bundle $E$
$$
P(E;\tau,z)=
\exp\biggl(-\sum_{n\ge 2}\frac{\wp^{(n-2)}(\tau,z)}{(2\pi i)^n n!}
\bigl(\sum_{i=1}^r x_i^n\bigr) \biggr).
$$
The second factor is the Witten factor 
$$
W(M;\tau)=
\exp\biggl(2\sum_{k\ge 2}
\frac{G_{2k}(\tau)}{(2k)!}
\bigl(\sum_{j=1}^d {x'_j}^{2k}\bigr) \biggr)
$$
which determines the Witten genus of the manifold $M$
as a function in one variable $\tau$ (see \S 1.3).
The derivation  of  order $(n-2)$ of the Weierstrass $\wp$-function
is  a meromorphic Jacobi form of weight $n$ and index $0$ with pole
of order $n$ along $z=0$.
Thus the coefficient of a monomial in $x_i$, $x'_j$ of the total 
degree $d$ in \thetag{1.8} is a meromorphic Jacobi form 
of weight $0$ and index $r/{2}$ with pole of order not bigger
than $(d-r)$. Therefore the product
$\eta(\tau)^{r-d}\vth(\tau,z)^{d-r}\chi(M,E,\tau,z)$ is holomorphic 
on $\bh_1\times \bc$ and it has trivial character.
It is a  weak Jacobi form since
its Fourier expansion
does not contain negative powers of $q$. 
Theorem 1.2 is proved.

We note that the special cocycle
defined by the logarithmic derivative of the Jacobi theta-series
is very natural in the following context.

\proclaim{Proposition 1.5}
Let $\phi\in J_{k,m}^{mer}$ be a (meromorphic) Jacobi form. 
The following functions
$$
\Psi(\tau, z)=\exp\bigl(-8\pi^2m\, G_2(\tau)z^2 \bigr)
\phi(\tau, z)
$$
and
$$
\Phi(\tau,z,\xi)=\exp\biggl(-8\pi^2m\, G_2(\tau)\xi^2
-\frac{\phi_z(\tau,z)}{\phi(\tau,z)}\xi\biggr)
\phi(\tau, z+\xi).
\tag{1.9}
$$
where $\phi_z=\pd{\psi}{z}$,
transform like a  modular form in $\tau$ and a Jacobi form of index $m$
in $\tau$ and $z$ respectively. More exactly they satisfy the equations
$$
\align
\Psi(\frac{a\tau+b}{c\tau+d},\,\frac z{c\tau+d})&=
(c\tau+d)^k
\Psi(\tau,z)
\\
\vspace{1\jot}
\Phi
(\frac{a\tau+b}{c\tau+d},\,\frac z{c\tau+d}, \frac \xi{c\tau+d})&=
v(\gamma)(c\tau+d)^k e^{2\pi i m 
\tsize{\frac { c z^2}{c\tau+d}}}
\Phi(\tau,z,\xi)\\
\vspace{1\jot}
\Phi(\tau, z+\lambda \tau+ \mu, \xi)&=
(-1)^{(\lambda+\nu)}
e^{-2\pi i m(\lambda^2 \tau+2\lambda z)}\Phi(\tau,z,\xi)
\endalign
$$
for any 
$\gamma=\left(\smallmatrix a&b\\c&d\endsmallmatrix\right) \in SL_2(\bz)$
and $\lambda, \mu\in \bz$.
\endproclaim

As a direct corollary we obtain that the coefficients
$g_n(\tau)$ and  $f_n(\tau,z)$ 
in the Taylor expansion of $\Psi$ in $z$ and  $\Phi$ in $\xi$
$$
\Psi(\tau,z)=\sum_{n\in \bz} g_n(\tau) z^n, \qquad
\Phi(\tau,z,\xi)=\sum_{n\in \bz} f_n(\tau,z)\xi^n
\tag{1.10}
$$
are a modular forms of weight $n+k$
and  a (meromorphic) Jacobi form  of weight $k+n$ and index $m$
respectively.
Proposition 1.4 gives us the exact ``additional" formula
for $\vth(\tau, z+\xi)$.
\smallskip

{\bf Examples for dimension $d=2,3,4$ and $d=6$, $r=2$.}
For $d=2$ and $d=3$ we consider an arbitrary vector bundle $E$ over 
$M$.
It follows  from the representation  \thetag{1.8} that
$$
\chi(M_2, E_r; \tau,z)= 
\frac{1}{2(2\pi i)^2}\biggl(\sum_{i=1}^{r} x_i^2\biggr) [M]
\frac {\vth^r}{\eta^{r+4}}\wp=
\frac{1}{24}(c_1(E)^2-2c_2(E))\,\phi_{0,1}
\biggl(\frac {\vth}{\eta}\biggr)^{r-2}
$$
where 
$\phi_{0,1}(\tau,z)= (y+10+y^{-1})+q(\dots)$
is the unique, up to a constant, weak Jacobi form of weight
$0$ and index $1$ (see \S 2). Let 
$\phi_{0,\frac{3}2}(\tau,z)= {\vth(\tau,2z)}/{\vth(\tau,z)}$.
For $d=3$   we have
$$
\chi(M_3, E_r; \tau,z)=
\frac{1}6 \bigl(c_1(E)^3-3c_1(E)c_2(E)+3c_3(E)\bigr) 
\phi_{0,\frac{3}2}(\tau,z)
\biggl(\frac {\vth}{\eta}\biggr)^{r-3}.
$$
Let us consider the case $d=4$, $r\ge 4$ with  $p_1(E_4)=p_1(M_4)$
($c_1(E_4)$ is arbitrary).
We have the following formula for MWG
$$
\chi(M_{4}, E_r; \tau, z)=
\bigl[ \hat A(M_4)\psi_{0,2}^{(2)}(\tau,z)+
\bigl(r\hat A(M_4)-\hat A(M_4, E_r)\bigr)\phi_{0,2}(\tau,z)\bigr]
\biggl(\frac{\vth(\tau,z)}{\eta(\tau)}\biggl)^{r-4},
\tag{1.11}
$$
where  $\phi_{0,2}(\tau,z)=y+4+y^{-1}+q(\dots)$ is defined in
\thetag{2.1} and 
$\psi_{0,2}^{(2)}=\phi_{0,1}^2-24\phi_{0,2}$
(see Theorem 2.2).

\smallskip
Let us consider the case $d=6$, $r=2$, $p_1(E)=p_1(T_M)$. We assume
also that MWG is holomorphic.
In the case under consideration this condition is equivalent 
to the relation
$<c_1(E)^4c_2(E)-\frac 1{6}c_1(E)^6>[M]=0$. Then
$$
\chi(M_6, E_2;\tau,z)=
\bigl[\hat A(M)E_6(\tau)\phi_{-2,1}+(2\hat A(M)- 
\hat A(M,E))E_{4,1}(\tau,z)\bigr]\eta(\tau)^{-8},
\tag{1.12}
$$
where $\phi_{-2,1}=\vth^2/\eta^6$, $E_6(\tau)$ is the Eisenstein series
of weight $6$ and $E_{4,1}(\tau,z)$ is the Eisenstein--Jacobi
series of weight $4$ and index $1$. ($E_{4,1}$ is equal to the Jacobi
theta-series of the lattice $\Bbb E_8$.) 
For more details and examples see \cite{G2}.

Theorem 1.2 shows us that  some possible ``{\bf strange periodicity}" of MWG
could exist. We note that 
if $p_1(E)=p_1(M)$ and $c_1(E)=0$, then MWG is the partition
function of $(0,2)$-symmetric non-linear sigma-model (see \cite{KM}).

\proclaim{Corollary 1.6}
Let us assume  that 
$
\chi(M_1, E_2, \tau,z)\equiv \chi(M_1, E_2, \tau,z)
$.
Then 
$r_1=r_2$ and  $d_1\equiv d_2\mod 12$.
\endproclaim
(See also Question 2 and 3 after Theorem 2.3.)

\head
\S 2. $\bz$-structure of the graded ring of  Jacobi forms
and the special values of the elliptic genus
\endhead

The  structure over $\bc$ of the bigraded ring of all weak Jacobi forms 
was determined in \cite{EZ}.
The elliptic genus of a Calabi--Yau manifold is a weak Jacobi
form of weight $0$ with integral Fourier coefficients.
Thus one can put a question about the $\bz$-structure of 
the graded ring
$$
J^\bz_{0,*}=\bigoplus_{m\in \bz}J^\bz_{0,m}
$$
of all Jacobi forms  of weight $0$ of integral weight
with integral Fourier coefficients.
For an arbitrary vector bundle $E$, the form 
$\chi(M,E; \tau,z)$ can be written, up to a $\eta$-factor,
 as a linear combination of weak Jacobi forms of positive weight
(see \thetag{1.12}). Thus we have also a question about 
the $\bz$-structure of the bigraded ring of all weak Jacobi
forms of arbitrary weight and index with integral Fourier coefficients. 
We introduce an ideal of $J^\bz_{0,*}$
$$
J^\bz_{0,*}(q)=\{\phi \in J^\bz_{0,*}\,|\, 
\phi(\tau,z)=\sum_{n\ge 1, \ l\in \bz} a(n,l) q^n y^l\}
$$
consisting of  the  Jacobi forms without $q^0$-term.
Using standard considerations with divisors of one can prove

\proclaim{Lemma 2.1}  {\bf A}.
Let  $m$ be integral, then we have 
$$
J_{2k,m+\frac 1{2}}^\bz=\phi_{0, \frac 3{2}}\cdot J_{2k,m-1}^\bz,
\qquad
J_{2k+1,m+\frac 1{2}}^\bz=\phi_{-1, \frac 1{2}}\cdot J_{2k+2,m}^\bz
$$
where 
$\phi_{0, \frac 3{2}}(\tau, z)={\vth(\tau,2z)}/{\vth(\tau,z)}$ and 
$\phi_{-1, \frac 1{2}}={\vth(\tau,z)}/{\eta(\tau)^3}$.

{\bf B}.
The ideal $J^\bz_{0,*}(q)$ is  principal.
It is  generated by a weak Jacobi form of weight $0$  and index $6$
$$
\xi_{0,6}(\tau,z)=\Delta(\tau)\phi_{-1,\frac {1}2}(\tau,z)^{12}=
\frac{\vartheta(\tau,z)^{12}}{\eta(\tau)^{12}}=
q(y^{\frac{1}2}-y^{-\frac{1}2})^{12}+q^2(\dots).
$$
\endproclaim
There exists only one (up to a constant) weak Jacobi form of 
weight $0$ and index $1$ 
$$
\phi_{0,1}(\tau,z)=-\frac{3}{\pi^2}\frac {\wp(\tau,z)\vth(\tau,z)^2}
{\eta(\tau)^6}=
(y+10+y^{-1})+q(10y^{\pm 2}-88y^{\pm 1}-132)+\dots
$$
(see \cite{EZ}).
In the theory of generalized Lorentzian Kac--Moody algebras
(see \cite{GN1--GN4}) we defined  the following important
Jacobi forms of small index:
$$
\align
\phi_{0,2}(\tau ,z)&=
{\tsize\frac{1}2} \eta(\tau )^{-4}
\sum_{m\,,n\in \bz}
{(3m-n)}\biggl(\frac{-4}m\biggl)\biggl(\frac{12}n\biggl)
q^{\frac{3m^2+n^2}{24}}y^{\frac{m+n}2}\\
{}&=(y+4+y^{-1})+q(y^{\pm 3}-8y^{\pm 2}-y^{\pm 1}+16)
+q^2(\dots),
\tag{2.1}\\
\vspace{1\jot}
\phi_{0,3}(\tau ,z)&=\phi_{0,\frac{3}2}^2(\tau ,z)=
(y+2+y^{-1})+q(-2y^{\pm 3}-2y^{\pm 2}+2y^{\pm 1}+4)+q^2(\dots),\\
\vspace{1\jot}
\phi_{0,4}(\tau ,z)&=
\frac {\vartheta(\tau ,3z)}{\vartheta(\tau ,z)}
=(y+1+y^{-1})-q(y^{\pm 4}+y^{\pm 3}-y^{\pm 1}-2)
+q^2(\dots).
\tag{2.2}
\endalign
$$
We have the following interpretation of these forms
as elliptic genus
$$
2\phi_{0,1}(\tau,z)=EG(K3, \tau,z),
\qquad 
\phi_{0,\frac{3}2}=EG(CY_3^{(e=2)}, \tau,z)
$$
($CY_3^{(e=2)}$ is  a Calabi--Yau $3$-fold with Euler number equals $2$.)
The form $\xi_{0,6}$ is the elliptic genus
of the trivial vector bundle of rank $12$ over a point
$$
\xi_{0,6}(\tau,z)=\chi(\cdot, \bc^{12};\tau,z).
$$
One can also  represent  these  functions as  symmetric polynomials  
in  the quotients of the Jacobi  theta-series  
$\vartheta_{ab}(\tau,z)/\vartheta_{ab}(\tau,0)$ of level $2$.
Let us put
$$
\xi_{00}=
\frac{\vartheta_{00}(\tau,z)}{\vartheta_{00}(\tau,0)},
\quad 
\xi_{10}=\frac{\vartheta_{10}(\tau,z)}{\vartheta_{10}(\tau,0)},
\quad
\xi_{01}=\frac{\vartheta_{01}(\tau,z)}{\vartheta_{01}(\tau,0)}.
$$
Then we have
$$
\gather
\phi_{0,1}(\tau,z)=4(\xi_{00}^2+\xi_{10}^2+\xi_{01}^2),
\qquad
\phi_{0,\frac 3{2}}(\tau,z)=
4 \xi_{00}\xi_{10}\xi_{01}
\\
\vspace{2\jot}
\phi_{0,2}(\tau,z)=2\bigl((\xi_{00}\xi_{10})^2+
(\xi_{00}\xi_{01})^2+(\xi_{10}\xi_{01})^2\bigl).
\endgather
$$
(To  check these formulae one should compare  only
$q^0$-terms of corresponding Jacobi forms.)

In the next theorem we construct a basis of   the module 
$J^\bz_{0,m}/J^\bz_{0,m}(q)$
and we find generators of the graded ring $J_{0,*}$.
\proclaim{Theorem 2.2} {\bf 1}. Let $m$ be  a positive integer.
The module
$$
J^\bz_{0,m}/ J^\bz_{0,m}(q)=
\bz\,[\psi_{0,m}^{(1)},\dots ,\psi_{0,m}^{(m)}]
$$
is a free $\bz$-module of rank $m$. Moreover
we can chose a basis  with the following $q^0$-terms 
$$
\align
[\psi_{0,m}^{(n)}(\tau,z)]_{q^0}&
=y^n-n^2 y+(2n^2-2)- n^2 y^{-1}+ y^{-n}\qquad\qquad
(2\le n\le m),\\
[\psi_{0,m}^{(1)}]_{q^0}&=\frac{1}{(12,m)}\,
\bigl(my+(12-2m)+my^{-1}\bigr)
\endalign
$$
where $(12,m)$ is the greatest common divisor of $12$ and $m$.

{\bf 2}. The graded ring of all  weak Jacobi forms
of weight $0$ with integral Fourier coefficients is finitely generated
$$
J^\bz_{0,*}=
\bigoplus_{m}J^\bz_{0,m}=
\bz\,[\phi_{0,1}, \phi_{0, 2}, \phi_{0,3}, \phi_{0,4}]
$$
where  $\phi_{0,1}$, $\phi_{0, 2}$, $\phi_{0,3}$
are algebraicly independent and 
$$
4\phi_{0,4}=\phi_{0,1}\phi_{0, {3}}-\phi_{0,2}^2.
$$
\endproclaim

The second claim of the theorem  is a corollary of the first part
which on can  prove by induction on $m$ and $n$.
We give here only the formulae for the most important exceptional 
Jacobi forms having the $q^0$-term of  type
$y+c+y^{-1}$ (see a detailed proof of Theorem 2.2 in \cite{G1}):
$$
\align
\phi_{0,6}(\tau,z)&=\phi_{0,2}\phi_{0,4}-\phi_{0,3}^2=
(y+y^{-1})+q(\dots),\\
\phi_{0,8}(\tau,z)&=\phi_{0,2}\phi_{0,6}-\phi_{0,4}^2
=(2y-1+2y^{-1})+q(\dots),\\
\phi_{0,12}(\tau,z)&=\phi_{0,4}\phi_{0,8}-2\phi_{0,6}
^2=(y-1+y^{-1})+q(\dots).
\endalign
$$
We note also that 
$$
\xi_{0,6}=-\phi_{0,1}^2\phi_{0,4}
+9\phi_{0,1}\phi_{0,2}\phi_{0,3}-8\phi_{0,2}^3-27\phi_{0,3}^2.
\tag{2.3}
$$
To prove that $\phi_{0,1}$, $\phi_{0,2}$ and $\phi_{0,3}$ are algebraicly
independent one has to consider  values at $z=\frac{1}{2}$.
We have
$$
\phi_{0,2}(\tau, \frac 1{2})\equiv 2,\quad
\phi_{0,3}(\tau, \frac 1{2})\equiv 0, \quad
\phi_{0,4}(\tau, \frac 1{2})\equiv -1.
$$
(The two  last identities follow from definition and the first
one is a corollary of the torsion relation of the theorem.) 
The restriction of 
$$
\phi_{0,1}(\tau, \frac 1{2})
=\alpha(\tau)=8+2^8q+2^{11}q^2+11\cdot 2^{10}q^3+
3\cdot 2^{14}q^4+359\cdot 2^9 q^5+\dots
$$
is a modular function with respect to $\Gamma_0(2)$ with  a character
of order $2$.

\example{Generators $\phi_{0,1}, \dots, \phi_{0,4}$ 
and Tate curves}
One can associate with  a plane cubic model of an  elliptic curve 
the quantities $b_2, b_4, b_6, b_8$ (see \cite{T}).
Theorem 2.2 shows us that the generators of the graded ring of the weak
Jacobi forms are  analytic analogies of the Tate's parameters since
$4b_8=b_2b_6-b_4^2$. The formula \thetag{2.3} is the expression for
the discriminant $\Delta$ of the Tate curve in terms of $b_i$.
\endexample

\example{Minimal models}As we mentioned above the elliptic genus
appeared in physics as the partition function of 
a $N=2$ superconformal field theory (SCFT). Among them we have 
the elliptic genus a sigma model with  a Calabi--Yau manifold
as  targent space, the elliptic genus of Landau--Ginzburg models
and of $N=2$ minimal models (see \cite{W2}, \cite{FY}, \cite{KYY}).

Let us consider   SCFT with central charge $\hat c$
($\hat c=c/3$). Then its elliptic genus
$Z(\tau,z)$  satisfies the functional equations
$$
Z\biggl(\frac{a\tau+b}{c\tau+d},\,\frac z{c\tau+d}\biggr)=
\,e^{\pi i \hat c\tsize{\frac { c z^2}{c\tau+d}}}
\,Z(\tau,z)
\qquad (
\gamma=\left(\smallmatrix a&b\\c&d\endpmatrix \in SL_2(\bz))
$$
and 
$$
Z(\tau, z+\lambda \tau+ \mu)=
(-1)^{2\hat c(\lambda+\mu)}
\,e^{-2\pi i \hat c (\lambda^2 \tau+2\lambda z)}Z(\tau,z)
\qquad (\lambda, \mu\in h\bz) 
$$
where $h\hat c$ is integer. 
(One can determine $h$ in terms of the $U(1)$ 
charges of the chiral ring elements in some sector.) 
We would like to note that the transformation 
$z\to hz$ is a Hecke type operator which multiplies the index
by $h^2$:
$$
\phi(\tau,z)\in J_{k,m}\to \phi(\tau,hz)\in J_{k, h^2m}.
$$
(See \cite{GN4, \S 1.4} where such operators were used
in construction of  Siegel theta-series.)
Therefore 
$$
Z(\tau, hz)\in J_{0, h^2 \hat c}.
\tag{2.4}
$$
Let us consider the elliptic genus 
$Z^{(min)}_k(\tau,z)$ of the $N=2$ minimal model with 
$\hat c=\frac{k}{k+2}$, where $k=1,2,\dots$.
Thus 
$$
Z^{(min)}_k(\tau, (k+2)z)\in J^{w, \bz}_{0,\frac 1{2} k(k+2)}.
$$
In particular, for $k=1$ the elliptic genus of $N=2$ minimal model
is equal, up to a constant, to the Jacobi form $\phi_{0,3/2}$ 
which is the elliptic genus of a Calabi--Yau $3$-fold.
If $k=2$, then the elliptic genus is a linear combination
of the four  forms $\phi_{0,4}$, $\psi_{0,4}^{(2)}$,
$\psi_{0,4}^{(3)}$  and $\psi_{0,4}^{(4)}$ constructed in Theorem 2.2.
The corresponding coefficients are determined by $q^0$-term
of the elliptic genus.

For $k=3$ the index is equal to $\frac {15}{2}$. Thus the elliptic genus
is a combination of the seven forms
$\phi_{0,3/2}\cdot \psi_{0,6}^{(m)}$ ($m=1,\dots,6$)
and $\phi_{0,3/2}\cdot\xi_{0,6}$.
This is the first case when one needs an information about one coefficient
of the $q^1$-term in order to find the elliptic genus.
We note that if $k$ is odd, then $Z^{(min)}_k(\tau, (k+2)z)$
is divisible by $\phi_{0,3/2}$.
\endexample

\example{Cobordism ring} G. H\"ohn was calculated 
the cobordism ring $\Omega^{SU}_*\otimes \bq$ over $\bq$
of $SU$-manifolds (see  \cite{H\"o}).
Using Theorem 2.2 above we can solve this problem over $\bz$.
\endexample

We have  also a result about the structure of the bigraded ring
of all integral weak Jacobi forms 
$$
J^{w,\bz}_{*,*/2}
=\bigoplus_{k\in \bz, \,m\in \frac{1}{2}\bz}J^{w, \bz}_{k,m}.
$$
\proclaim{Theorem  2.3}
$$
J_{*, */2}^{w, \bz}=\bz[E_4, E_6, \Delta,
E_{4,1}, E_{4,2}, E_{4,3}, 
E_{6,1}, E_{6,2}, E'_{6,3}, 
\phi_{0,1}, \phi_{0,2}, \phi_{0,\frac{3}{2}}, \phi_{0,4},
\phi_{-1, \frac 1{2}}]
$$
where 
$E_4(\tau)$, $E_6(\tau)$ and $\Delta(\tau)$ are generators of the graded 
ring of $SL_2(\bz)$-modular forms,
$\phi_{-2, 1}=\vth^2/\eta^6$,
$E_{4,1},\dots E_{6,3}$ are the Eisenstein--Jacobi series with the
zeroth Fourier coefficient equals to $1$ 
and 
$E'_{6,3}=E_{6,3}+\frac{22}{61}\Delta_{12}\phi_{-2,1}^3$.
\endproclaim

As we mentioned at the end of \S 1, the MWG is a partition function
of $(0,2)$-symmetric sigma model.
We would like to formulate  a problem  about existence of
some special vector bundles with an  interesting  $(0,2)$-partition function.
We assume that  $M_d$ is a  complex compact  manifold 
of  complex dimension $d$, $E_r$ is a holomorphic  vector bundle of rank $r$
over $M$. Moreover 
$$
p_1(E_r)=p_1(M_d) \quad\text{ and }\quad c_1(E_r)= 0 \ (\text{over } \br).
$$
{\bf Question 1}. {\it The Eisenstein type.}
Let $d=6$, $r=2$ and
$$
\hat A(M_6)=0,      \qquad  \hat A(M_6, E_2)=-2.
$$
Similar questions  for dimensions $8$ and $10$:
$$
\hat A(M_8)=\hat A(M_8, E_4)=0,  \quad 
\hat A(M_8, \wedge^2E_4)=-2 
$$
and 
$$
\hat A(M_{10})=\hat A(M_{10}, E_6)=\hat A(M_{10}, \wedge^2E_6)=0,\ 
\hat A(M_{10}, \wedge^3E_6)=-2.
$$
The MWG of a  vector bundle with invariant as above is 
an  Eisenstein-Jacobi series from Theorem 2.3.

{\bf Question 2}. {\it The ``trivial" type.}
Let $d=12$, $r=2$ and
$$
\hat A(M_{12})=\hat A(M_{12}, E_2)=0, \qquad 
\hat A(M_{12}, E_2\otimes T_M)\ne 0.
$$

For a $K3$-surface we can formulate

{\bf Question 3}. {\it The K3 type.}
Let  $d=14$, $r=2$,
$$
\hat A(M)=\hat A(M, E)=0, \qquad \hat A(M, E\otimes T_M)=-10.
$$
The last condition is equivalent to 
$\hat A(M, S^2E)+ 2\hat A(M, T_M)=2$.

\smallskip
{\bf Question 4}. {\it The Kac--Moody type.}
Let $d=14$, $r=2$.
$$
\hat A(M)=0, \quad \hat A(M, E)=-2
$$
and 
$$
\hat A(M, E\otimes T_M)=-68.
$$
The last condition is equivalent to 
$
\hat A(M, S^2E)+ 2\hat A(M, T_M)=0
$.
Such a vector bundle (if it exist) will be related to the simplest
hyperbolic Kac--Moody algebra.
(See \S 4 bellow.)

\medskip

Using Theorems 2.2 and 2.3   we can  analyze the value  
of the elliptic genus 
at the following special points
$z=0$ (Euler number), $z=\frac 1{2}$ (signature),
$z=\frac{\tau+1}2$ ($\hat A$-genus) and $z=\frac 1{3},
\ \frac 1{4},\ \frac 1{6}$.
For this end we have to study the restriction of the 
generators of the graded ring of the integral week Jacobi forms.
A special value of a Jacobi form
is a modular form  in $\tau$.
In the next lemma we give a little more
precise statement than in  \cite{EZ, Theorem 1.3}.
\proclaim{Lemma 2.4} Let $\phi\in J_{0,t}$ ($t\in \bz/2$) and 
$X=(\lambda, \mu)\in \bq^2$. Then
$$
\phi|_{X}(\tau,0)=\phi(\tau, \lambda\tau+\mu)
\exp(2\pi i t(\lambda^2\tau+ \lambda \mu))
$$ 
is an automorphic form of weight $0$ with a character 
with respect to  the  subgroup 
$$
\Gamma_X=\{M\in SL_2(\bz)\,|\, XM-X\in \bz^2\}.
$$ 
\endproclaim

It is easy to see that if $\phi\in J_{k,m}^\bz$ with integral $m$, then 
the form $\phi(\tau,\frac 1{N})$ still has integral Fourier coefficients if
$N=1, \dots, 6$. In particular, the value of $\xi_6(\tau,z)$ at 
these points is related to the ``Hauptmodul" for the corresponding
group $\Gamma_0(N)$:
$$
\aligned
\xi_6(\tau,\frac{1}2)&=2^{12}\frac{\Delta(2\tau)}{\Delta(\tau)},\\
\xi_6(\tau,\frac{1}3)&=
3^{6}\left(\frac{\Delta(3\tau)}{\Delta(\tau)}\right)^{1/2},
\endaligned
\qquad
\aligned
\xi_6(\tau,\frac{1}4)&=
2^{6}\left(\frac{\Delta(4\tau)}{\Delta(\tau)}\right)^{1/2},
\\
\xi_6(\tau,\frac{1}6)&=
\left(\frac{\Delta(\tau)\Delta(6\tau)}
{\Delta(2\tau)\Delta(3\tau)}\right)^{1/2}.
\endaligned
$$
Let us analyze the corresponding values of the four generators $\phi_{0,n}$
of the graded ring $J_{0,*}^{\bz}$. From the definition 
(see \thetag{2.1}--\thetag{2.2}) and the identity
$4\phi_{0,4}=\phi_{0,1}\phi_{0, {3}}-\phi_{0,2}^2$
we obtain 
$$
\phi_{0,1}(\tau, 0)=12,\quad 
\phi_{0,2}(\tau, 0)=6,\quad
\phi_{0,3}(\tau, 0)=4,\quad
\phi_{0,4}(\tau, 0)=3
\tag{2.5}
$$
and
$$
\aligned
\phi_{0,1}(\tau, \frac{1}2)&=\alpha(\tau)\\
\phi_{0,2}(\tau, \frac{1}2)&=2\\
\phi_{0,3}(\tau, \frac{1}2)&=0\\
\phi_{0,4}(\tau, \frac{1}2)&=-1
\endaligned
\qquad
\aligned
\phi_{0,1}(\tau, \frac{1}3)&=\beta^2(\tau)\\
\phi_{0,2}(\tau, \frac{1}3)&=\beta(\tau)\\
\phi_{0,3}(\tau, \frac{1}3)&=1\\
\phi_{0,4}(\tau, \frac{1}3)&=0
\endaligned
\qquad
\aligned
\phi_{0,1}(\tau, \frac{1}4)&=\frac{\gamma(\tau)^4+4}{\gamma(\tau)}\\
\phi_{0,2}(\tau, \frac{1}4)&=4\gamma^2(\tau)\\
\phi_{0,3}(\tau, \frac{1}4)&=2\gamma(\tau)\\
\phi_{0,4}(\tau, \frac{1}4)&=1.
\endaligned
\tag{2.6}
$$
The automorphic functions $\alpha(\tau)$, $\beta(\tau)$ and $\gamma(\tau)$
are automorphic forms of weight $0$ with respect to the group
$\Gamma_0$, $\Gamma_0^{(1)}(3)$ and $\Gamma_0^{(1)}(4)$ respectively.
These functions have integral Fourier coefficients.
The identity \thetag{2.3} gives us the following relations
between  the automorphic functions  $\alpha$, $\beta$ and $\gamma$
$$
\gather
2^{12}\frac{\Delta(2\tau)}{\Delta(\tau)}=
\alpha(\tau)^2-64, \qquad
3^{6}\left(\frac{\Delta(3\tau)}{\Delta(\tau)}\right)^{1/2}
=\beta(\tau)^3-27\\
2^{6}\left(\frac{\Delta(4\tau)}{\Delta(\tau)}\right)^{1/2}
=4(\left(\frac{\gamma(\tau)}2\right)^2
-\left(\frac 2{\gamma(\tau)}\right)^2).
\endgather
$$
It follows that
$$
\alpha(\tau)-8\equiv 0 \mod 2^8,\qquad
\beta(\tau)-3\equiv 0 \mod 3^3
\tag{2.7}
$$
(compare with the formula for $\phi_{0,1}(\tau,\frac 1{2})$ above).
Using the  definition of $\phi_{0,3}$ and $\gamma(\tau)$ 
and the relations between 
the Jacobi theta-series $\vartheta_{ab}$ of level $2$  we have 
$$
\gamma(\tau)=
\frac{\vartheta_{00}(2\tau)}{\vartheta_{01}(2\tau)}=
\frac{\vartheta_{00}(2\tau,0)}{\vartheta_{01}(2\tau,0)}.
$$
One can  check that
$
\phi_{0,1}(\tau, 2z)=\phi_{0,2}^2(\tau,z)-8\phi_{0,4}(\tau,z)
$.
Thus
$$
\alpha(\tau)=16\gamma(\tau)^4-8=
16\frac{\vth_{00}^4(2\tau)}{\vth_{01}^4(2\tau)}-8.
$$
In particular {\it all Fourier coefficients of $\gamma(\tau)$
and $\alpha(\tau)$ are positive}.

\example{Example 2.5}{\bf $\hat A$-genus.}
Let $X=(\frac{1}{N}, \frac{1}{N})$.
Then  $\Gamma_X$ (see Lemma 2.4) contains the principle congruence subgroup
$\Gamma_1(N)$. In some cases $\Gamma_X$ will be strictly  larger.
For example, if $X_2=(\frac{1}{2}, \frac{1}{2})$,  then
$$
\phi|_{X_2}(\tau,0)=\phi(\tau,\frac{\tau+1}2)
\exp(\frac{\pi i}2({\tau+1}))
$$
is an automorphic form with respect of the so-called theta-group
$$
\Gamma_\theta=
\biggl\{ M=\pmatrix a&b\\c&d\endpmatrix \in 
 SL_2(\bz)\,|\  M\equiv \pmatrix 1&0\\0&1\endpmatrix 
\ \text{ or }\ \pmatrix 0&1\\1&0\endpmatrix
\mod 2
\biggr\}.
$$
The corresponding character is given by
$\epsilon_2(M)=\exp(2\pi i m(d+b-a-c)/4)=\pm 1$. 
This  character is trivial if index $m$ of Jacobi
form is even. 
Let us  consider $\Gamma_\theta$-automorphic  function
$$
\hat \phi_{m}(\tau)=
q^{-\frac{m}4}\phi_{0,m}(\tau, -\frac{\tau+1}2).
$$
We have
$$
\hat\phi_3= 0,\quad \hat\phi_4= -1,\quad
\hat\phi_2= -2,\quad
\hat\xi_6=\hat \phi_1^2+64=\left(\frac{\vth_{00}}{\eta}\right)^{12}
$$
where  
$$
\hat\phi_1(\tau)=4\frac{\vth_{10}^4-\vth_{01}^4}{\vth_{01}^2\vth_{10}^2}
=-q^{-\frac 1{4}}+20q^{\frac 1{4}}+\dots 
\in \frak M_0^\bz(\Gamma_{0}(2), \epsilon_2).
$$
\endexample

\medskip
Now we  analyze some special values  of the elliptic genus.
As it easy follows from \thetag{1.2} we get Euler number of a 
Calabi--Yau manifold $M_d$ for $z=0$ ($d$ is arbitrary) and 
 and its signature for $z=\frac{1}2$ 
($d$ is even):
$$
\align
\chi(M_d, \tau, 0)&=e(M_d),\\
\chi(M_d, \tau, \frac{1}2)&=\sigma_M(\tau)=
(-1)^{\frac{d}2}s(M_d)+q(\dots)\in 
\frak M_0^\bz(\Gamma_{0}(2), v_2),
\quad
v_2(\left(\smallmatrix a&b\\c&d\endpmatrix)=
e^{\pi i m \frac {c}2}.
\endalign
$$
The formulae \thetag{2.5} gives  us some divisibility
of Euler number of Calabi--Yau manifolds.
We note that the quotient $e(M)/24$ appears in physics as obstruction
to cancelling the tadpole (see \cite{SVW} where it was proved
that $e(M_4)\equiv 0 \mod 6$).
\proclaim{Proposition 2.6}Let $M_d$ be an almost complex manifold 
of complex 
dimension $d$ such that $c_1(M)=0$ in $H^2(M, \br)$.
Then 
$$
d\cdot e(M_d)\equiv 0 \mod 24.
$$
If $c_1(M)=0$ in $H^2(M, \bz)$,  then we have a more strong congruence
$$
e(M)\equiv 0 \mod 8 \qquad\text{if }\  d\equiv 2 \mod 8.
$$
\endproclaim
\remark{Remark} More generally, we have the following relation 
for $\hat A$-genus of a spin bundle. 
Let us consider the case
$$
d=r,\quad p_1(E)=p_1(T_M), \quad c_1(E)=0\  (\text{over }\bz).
$$
We consider  spin bundle $\Delta^+(E)-\Delta^-(E)$ such that
$
\hbox{ch}(\Delta^+(E)-\Delta^-(E))
=\prod_{i=1}^r(e^{x_i/2}-e^{-x_i/2})
$. Then
$$
\hat A(M, \Delta^+(E)-\Delta^-(E))\equiv
\cases 0\mod 8& d\equiv 1,\,2,\,5 \mod 8\\
 0\mod 4& d\equiv \ \ \ 6\,,\,7 \mod 8\\
 0\mod 2& d\equiv \ \ \ \,3,\,4 \mod 8\\
 0\mod 3& d\not \equiv \ \ \ \ \ \ \ 0 \mod 3.
\endcases
$$
\endremark
\demo{Proof}The first fact follows simply from  \thetag{2.5}.
If $d\equiv 2 \mod 8$
one can write the elliptic genus as a polynomial  over $\bz$ in 
the generators $\phi$
$$
e(M_d)\equiv P(\phi_{0,1}, \phi_{0,2}, \phi_{0,3}, \phi_{0,4})|_{z=0}
\equiv c_{1,m}(\phi_{0,1}|_{z=0})(\phi_{0,4}|_{z=0})^{\frac {d-2}8}
\mod 8.
$$
If one put $z=-\frac{\tau+1}2$, i.e.,
$y=-q^{1/2}$, then one can see that the series 
$$
\Bbb E_{q,-q^{1/2}} = \bigotimes_{n\geq 1}
{\bigwedge}_{q^{n/2}}T_M
\otimes 
\bigotimes_{n\geq 1}{\bigwedge}_{q^{n/2}} {T}_M^*
\otimes 
\bigotimes_{n\geq 1}
S_{q^n} (T_M \oplus {T}_M^*)
$$
is $*$-symmetric.  According to the  Serre duality all Fourier 
coefficients of  $\hat\chi(M_d,\tau)$ are even.
The constant $ c_{1,m}$ from the last congruence
is equal to the coefficient of $\hat\chi(M_d,\tau)$ at 
the minimal negative power of $q$. Therefore
$c_{1,m}$ is even and we obtain divisibility of $e(M_{8m+2})$ by 8.
\enddemo
We note that divisibility of $d \cdot e(M)$ by $3$ was proved by 
F. Hirzebruch in 1960.
For  a hyper-K\"ahler compact  manifold the claim of the proposition above
was proved  by S. Salamon  in \cite{S}.
After my talk on the  elliptic genus  at a seminar of MPI
in Bonn in April 1997 Professor F. Hirzebruch informed me that 
the result of Proposition 2.6 was known for him (non-published).
Using some natural examples he also proved that 
this property of divisibility 
of the Euler number modulo $24$ is strict (see \cite{H2}).

Formulae \thetag{2.6} provide us with a formula for the signature 
$\chi(M_d;\tau, \frac{1}2)$ as a polynomial in $\alpha(\tau)$.
As a corollary of \thetag{2.6} and Theorem 2.2
we have that 
for an arbitrary Jacobi form of integral index
$$
\aligned
\phi_{0,4m}(\tau,\frac{1}2)&=c+2^{13}q(\dots)\\
\phi_{0,4m+2}(\tau,\frac{1}2)&=2 c+2^{12}q(\dots)
\endaligned
\qquad 
\aligned 
\phi_{0,4m+1}(\tau,\frac{1}2)&=8 c+2^{8}q(\dots)\\
\phi_{0,4m+3}(\tau,\frac{1}2)&=16 c+2^{9}q(\dots).
\endaligned
$$
Similar to the proof of Proposition 2.4 we obtain a better congruence
for the signature of a manifold with dim$\equiv 2 \mod 8$
and $c_1(M)=0$:
$$
\chi(M_{8m+2};\tau,z)=16c+2^{9}q(\dots). \tag{2.9}
$$
It  is interesting that the values of the Hirzebruch $y$-genus
at $y=e^{2\pi i/3}$ and $y=i$ also have some properties of divisibility.
For $z=\frac{1}3$ (resp. $z=\frac{1}4$) we can write 
$\phi_{0,m}(\tau, \frac{1}3)$ 
(resp. $\phi_{0,m}(\tau, \frac{1}4)$)
as a polynomial in $\beta(\tau)=3+27(q+\dots)$ 
(resp. in $\gamma(\tau)^{\pm 1}$).
This gives us the following results
$$
\gather 
\phi_{0,3m}(\tau,\frac{1}3)=c+3^{6}q(\dots),\qquad
\phi_{0,3m+1}(\tau,\frac{1}3)=9c+3^{4}q(\dots)\\
\phi_{0,3m+2}(\tau,\frac{1}3)=3c+3^{3}q(\dots).
\endgather
$$
Thus  we have 

\proclaim{Proposition 2.7} If $c_1(M)=0$ (over $\br$), then
$$
\gather
\chi(M_{6m}; \tau,\frac{1}3)\equiv c_1 \mod 3^6,\qquad
\chi(M_{6m+2}; \tau,\frac{1}3)\equiv 9c_2 \mod 3^4,\\
\chi(M_{6m+4}; \tau,\frac{1}3)\equiv 3c_3 \mod 3^3.
\endgather
$$
where $c_1, c_2, c_3\in \bz$. For $z=\frac 1{4}$ we have:
$$
\gather
\chi(M_{8m+2}; \tau,\frac{1}4)=4c+2^{4}q(\dots),\qquad
\phi_{0,4m+2}(\tau,\frac{1}4)=4c+2^{5}q(\dots)\\ 
\phi_{0,4m+3}(\tau,\frac{1}4)=2c+2^{8}q(\dots).
\endgather
$$
\endproclaim

\head
\S 3. $\hat A_2^{(2)}$-genus.
\endhead

The formal power series $\bold E_{q,y}$ over $K(M)$ (see \thetag{1.1})
is a geometric analog of the Jacobi theta-series $\vth(\tau,z)$
which is the Weyl--Kac denominator function of the affine Lie  algebra 
$\hat A_1^{(1)}$.
A similar construction we can propose for an arbitrary affine 
Lie algebra.
We are going  to  consider a general case in a separate publication.
In this section we consider the case of the affine Lie algebra 
$\hat A_2^{(2)}$.

Let us define
$$
\multline
{\bold E}_{q,y}^{(2)}=  {\bigwedge}_{y^{-1} } E^*
\otimes
\bigotimes_{n\ge 1}
{\bigwedge}_{-q^{n}y^{-2}}  \Psi_2(E^*) \otimes 
\bigotimes_{n \ge 1}{\bigwedge}_{-q^n y^2}  \Psi_2(E)\\
\otimes 
\bigotimes_{n\ge 1} S_{q^ny^{-1}} E^* \otimes 
\bigotimes_{n\ge 1} S_{q^ny} E \otimes 
\bigotimes_{n= 1}^{\infty} S_{q^n} (T_M \oplus T_M^*)
\endmultline
$$
where 
$\Psi_2(E)$ is the second Adams operation on vector bundle $E$.
We remind that 
$$
\hbox{ch}(\Psi_2(E))=\hbox{ch}(E\otimes E)-\hbox{ch}(\wedge^2 E)=
\sum_{i=1}^r e^{2x_i}.
$$
The series ${\bold E}_{q,y}$ is a geometric variant of the Jacobi triple
product and  ${\bold  E}_{q,y}^{(2)}$ 
is a geometric analog of the quintiple product
$$
\vartheta_{3/2}(\tau ,z)=
q^{\frac 1{24}}r^{-\frac 1{2}}
\prod_{n\ge 1}(1+q^{n-1}r)(1+q^{n}r^{-1})(1-q^{2n-1}r^2)
(1-q^{2n-1}r^{-2})(1-q^n).
$$
We note that
$$
\vartheta_{3/2}(\tau ,z)=\frac{\eta(\tau )\vartheta(\tau , 2z)}
{\vartheta(\tau , z)}\in 
J_{\frac 1{2}, \frac 3{2}}(v_\eta).
$$
For $\hat A_2^{(2)}$ we can give the following definition
(compare with Definition 1.1).
\definition{Definition 3.1}
{\it $ \hat A_2^{(2)}$-genus} 
of a complex  vector bundle $E$ of rank $r$
over a compact  complex manifold $M$ of dimension $d$
is defined as follows
$$
\align
{}&\alpha(M, E;\tau,z)=
q^{-d/12}y^{r/2}\int_M  
\exp \biggl(\frac 1{2}\bigl(c_1(E)-c_1(T_M)\bigr)\biggr)\\
\vspace{2\jot}
{}&
\exp\biggl(\bigl(3p_1(E)-p_1(T_M)\bigr)\cdot G_2(\tau)\biggr)
\exp\biggl(-\frac  {c_1(E)}{2\pi i} 
\dfrac{\partial\ }{\partial z}
\hbox{log} (\vth_{3/2}(\tau,z)) \biggr)
\,\hbox{ch}({\bold E}_{q,y}^{(2)})\, \hbox {td}(T_M).
\endalign
$$
\enddefinition

\proclaim{Theorem 3.1} Let $E$ be a holomorphic  vector bundle 
of rank $r$ over a compact complex manifold $M$ of dimension $d$. 
Let $\alpha (M,E;\tau, z)$ be  the $\hat A_2^{(2)}$-genus. 
Then 
$$
\eta(\tau)^{d+r}\vartheta_{3/2}(\tau ,z)^{d-r}
\alpha (M,E;\tau, z)\in J_{d, \frac 3{2}d}^{w}
$$
is a weak  Jacobi form of weight $d$ and index $\frac{3d}2$.
\endproclaim
The proof is similar to the proof of Theorem 1.2.

\head
\S 4. SQEG and hyperbolic root systems
\endhead

We can consider   $n$-fold symmetric product of the manifold  $M$,
i.e., the orbifold space $\Mn=M^n/S_n$, 
where $S_n$ is the symmetric group of $n$ elements.
This is a singular manifold but one can define the  string  orbifold
elliptic genus of $\Mn$ (see for details the talk 
of R. Dijkgraaf at ICM-1998 in Berlin \cite{D}).
Using some  arguments from the conformal field theory 
on orbifolds it was proved in \cite{DVV} and  \cite{DMVV} 
that the string elliptic genus
of the second quantization 
$\cup_{n\ge1} \Mn$ of  a Calabi--Yau manifold $M$
coincides with the second quantized elliptic genus of the 
given manifold:
$$
\sum_{n=0}^\infty p^n \chi_{orb}(\Mn;\, q,y)  
= \prod_{m\ge 0,\, l,\, n>0}
\frac 1{
(1- q^m y^l p^n)^{f(mn, l)}},
\tag{4.1}
$$
where 
$
\chi(M,\tau,z)=\sum_{m\ge 0,\ l\in \bz \ (or\,\bz/2)} f(m,l)q^my^l
$
is the elliptic genus of $M$.

For a $K3$ surface, the product in the left hand side of \thetag{4.1}
is essentially  the power $-2$  of the infinite product expansion 
of
the product of all even theta-constants
has found in \cite{GN1}.
Following \cite{DVV, \S 4}  we call the  product in \thetag{4.1}
 {\it second-quantized elliptic genus} (SQEG) of the manifold  $M$.
\proclaim{Theorem 4.1} (See \cite{GN4}, \cite{G1}.)
Let $M=M_d$ be a compact complex manifold
of dimension $d$
with trivial  $c_1(M)$,
$$
\chi(M;\, \tau,z)=\sum_{m\ge 0,\ l\in \bz \ (or\,\bz/2)} f(m,l)q^my^l
$$
be its elliptic genus and $\ \hbox{\rm SQEG}(M;\,Z)$  ($Z\in \bh_2$)
be its second quantized elliptic genus.
We define a factor
$$
H(M;\,Z)=
\cases
\eta(\tau)^{-\frac 1{2}(e-3\chi'_{d_0})}
\prod_{p=1}^{d_0}\bigl(\vartheta(\tau, pz)\,e^{\pi i p^2\omega}\bigl)
^{-\chi'_{d_0-p}}&\ \text{ if }\  d=2d_0\\
\eta(\tau)^{-\frac 1{2}{e}}
\prod_{p=1}^{d_0}\bigl(\vartheta(\tau, \frac{2p-1}2 z)\,
e^{\frac {1 }4 \pi i (2p-1)^2\omega}\bigl)
^{-\chi'_{d_0-p+1}}&\ \text{ if } \ d=2d_0+1
\endcases
$$
where $e=e(M)$ is Euler number of $M$ and 
$\chi'_p=(-1)^p\chi_p(M)$  (see \thetag{1.2}).
Then  the product
$$\align
E(M;\,Z)&=\Psi(M;\,Z)\cdot \hbox{\rm SQEG}(M;\,Z)\qquad (d=2d_0)\\
E^{(2)}(M;\,Z)&=(E |\Lambda_2) (M;\,Z)\qquad \qquad (d=2d_0+1)
\endalign
$$
determines a Siegel automorphic form of weight $-\frac 1{2}{\chi'_{d_0}(M)}$
if $d$ is even and of of weight $0$ if $d$ is odd 
with a character  or a multiplier system 
of order ${24}/(24, e)$ 
with respect to a double extension of the  paramodular group 
$\Gamma_d^+$ (resp. $\Gamma_{2d}^+$), if $d$ is even
(resp. $d$ is odd).
\endproclaim
\smallskip
{\bf The case of $\bold C\bold Y_{\bold 4}$.}
The basic  Jacobi modular forms for this dimension are 
the Jacobi forms $\phi_{0,2}$ and $\psi_{0,2}^{(2)}$
(see Theorem 2.2, part 1). They correspond to the following
cusp forms for the paramodular group  $\Gamma_2$ 
(see \cite{GN1} and \cite{GN4}):
$$
\align
\Delta_2(Z)&=\ml(\phi_{0,2}(\tau,z))=
\hbox{Lift}(\eta(\tau)^3\vartheta(\tau,z))\\
{}&=\sum_{N\ge 1}\
\sum\Sb
 n,\,m >0,\,l\in \bz\\
\vspace{0.5\jot} n,\,m\equiv 1\,mod\,4\\
\vspace{0.5\jot} 2nm-l^2=N^2
\endSb
\hskip-4pt
N\biggl(\frac {-4}{Nl}\biggr)
\sum_{a\,|\,(n,l,m)} \biggl(\frac {-4}{a}\biggr)
\, q^{n/4} y^{l/2} s^{m/2}
\in \frak M_2^{cusp}(\Gamma_2,v_\eta^6\times v_H)
\endalign
$$
and
$$
\Delta_{11}(Z)=\hbox{Lift}(\eta(\tau )^{21}\vartheta(\tau ,2z))=
\ml(\psi_{0,2}^{(2)}(\tau,z))
\in \frak N_{11}(\Gamma_{2}).
$$
For an arbitrary  Calabi--Yau $4$-fold $M_4$
we have the following formula for its
SQEG
$$
E(M_4;\,Z)=\Delta_{11}(Z)^{-\chi_0(M)} \Delta_{2}(Z)^{\chi_1(M)}.
\tag{3.2}
$$
We note that $\Delta_{2}(Z)^4$ is the first $\Gamma_2$-cusp form
with trivial character and $\Delta_{11}(Z)$ is the first cusp form of odd
weight with respect to $\Gamma_2$. 

The Fourier expansion of the cusp forms  
$\Delta_2(Z)$, $\Delta_{11}(Z)$ and 
$\frac {\Delta_{11}(Z)}{\Delta_{2}(Z)}$
coincide with  the Weyl--Kac--Borcherds denominator formula
of  generalized  Kac--Moody super-algebras with a system 
of simple real  roots of hyperbolic type determined by 
Cartan matrix $A_{1,II}$, 
$A_{2,II}$ and $A_{2,0}$ respectively:
$$
A_{2,II}=
\pmatrix\hphantom{-}{2}&{-2}&{-6}&{-2}\cr
{-2}&\hphantom{-}{2}&{-2}&{-6}\cr
{-6}&{-2}&\hphantom{-}{2}&{-2}\cr
{-2}&{-6}&{-2}&\hphantom{-}{2}\cr
\endpmatrix,
\quad
A_{2,0}=
\left(\smallmatrix\hphantom{-}2&-2&-2\\
-2&\hphantom{-}2&\hphantom{-}0\\
-2&\hphantom{-}0&\hphantom{-}2
\endpmatrix,
\quad
A_{2,I}=
\pmatrix\hphantom{-}{2}&{-2}&{-4}&\hphantom{-}{0}\cr
{-2}&\hphantom{-}{2}&\hphantom{-}{0}&{-4}\cr
{-4}&\hphantom{-}{0}&\hphantom{-}{2}&{-2}\cr
\hphantom{-}{0}&{-4}&{-2}&\hphantom{-}{2}\cr
\endpmatrix
$$
(see \cite{GN1}--\cite{GN4}).
Thus, the formula \thetag{3.2} gives us three particular cases
of Calabi--Yau $4$-folds of Kac--Moody type when 
the second quantized elliptic genus is a power of 
the denominator function of the corresponding  Lorentzian Kac--Moody
algebra:
$$
\align
E(M_4;\,Z)&=\Delta_{11}(Z)^{-\chi_0}
\qquad\qquad\text{if } \chi_1=0\\ 
E(M_4;\,Z)&=\biggl(\frac{\Delta_{11}(Z)}{\Delta_{2}(Z)}\biggr)^{-\chi_0} 
\qquad\text{if } \chi_0(M)=-\chi_1(M)\\
E(M_4;\,Z)&=\Delta_{2}(Z)^{\chi_1}\qquad\qquad\quad\text{if } \chi_0(M)=0.
\endalign 
$$
For more details and for the cases of $d>4$ see \cite{G1}.
In \cite{GN2--GN4} we started to classify all generalized 
Lorentzian Kac--Moody algebras whose denominator function
is an automorphic form on the Siegel upper-half plane of genus two. 

\proclaim{Conjecture}The Weyl--Kac--Borcherds denominator 
function of arbitrary automorphic Lorentzian Kac--Moody algebra
with a non-empty  system of real simple roots of rank $3$
is related to the second quantized elliptic genus of 
some holomorphic vector bundle.
\endproclaim

The Igusa modular form $\Delta_{35}(Z)$ determines the automorphic correction
of the simplest hyperbolic Kac--Moody algebra with the Cartan matrix
$$
A=\pmatrix\hphantom{-}2&\hphantom{-}0&-1\\
\hphantom{-}0&\hphantom{-}2&\hphantom-2\\
-1&-2&\hphantom{-}2
\endpmatrix.
$$
The form $\Delta_{35}(Z)$ would be  the SQEG of 
a vector bundle of rank $2$ over a manifold of dim$=14$ (if it would exist) 
with invariants given in Question 3 in \S 2.

\enddocument 
\end


\Refs 
\widestnumber\key{DMVV}

\ref
\key D 
\by R. Dijkgraaf
\paper The Mathematics of Fivebranes
\jour Documenta Mathem. ICM-1998
\yr 1998 
\endref 

\ref
\key DVV 
\by R. Dijkgraaf, E. Verlinde and H. Verlinde 
\paper Counting dyons in $N=4$ string theory 
\jour Nucl. Phys.
\vol  B484 
\yr 1997
\pages 543--561  
\endref 

\ref
\key DMVV 
\by R. Dijkgraaf, G. Moore,  E. Verlinde, H. Verlinde 
\paper Elliptic genera of symmetric products
and second quantized strings
\jour Commun. Math. Phys.
\vol 185
\yr 1997
\pages 197--209
\endref 


\ref
\key EZ
\by M. Eichler,  D. Zagier
\book The theory of Jacobi forms
\yr 1985
\publ Progress in Math. 55, Birkh\"auser
\endref


\ref
\key
EOTY
\by
 T. Eguchi, H. Ooguri, A. Taormina, S.-K. Yang,
\paper
Superconformal Algebras and String Compactification on Manifolds
with $SU(N)$ Holonomy
\jour
 Nucl. Phys. 
\yr 1989
\vol  B315
\pages 193--221
\endref

\ref
\key FY
\by P. Di Francesco, S. Yankielowicz 
\paper Ramond Sector Characters and 
$N=2$ Landau--Ginzburg models
\jour Nucl. Phys.
\vol B409
\yr 1993
\pages 1886--210
\endref

\ref
\key G1
\by V\. Gritsenko
\paper Elliptic genus of Calabi--Yau manifolds and Jacobi and Siegel
modular forms
\jour Algebra i Analiz (St. Petersburg Math. J.)
\vol 11:5
\yr 1999
\endref



\ref
\key G2
\by V\. Gritsenko
\paper Modified Witten genus
\inbook  Course of lectures at RIMS
\publ Kyoto
\yr 1998
\endref




\ref
\key GH
\by V. Gritsenko, K. Hulek
\paper Commutator coverings of Siegel threefolds
\jour  Duke Math. J.
\vol 94
\yr 1998
\pages 509--542
\endref


\ref
\key GN1
\by V.A. Gritsenko, V.V. Nikulin
\paper Siegel automorphic form correction of some Lorentzi\-an
Kac--Moody Lie algebras
\jour Amer. J. Math.
\yr 1997 
\vol 119
\pages  181--224
\endref



\ref
\key GN2 
\by V.A. Gritsenko, V.V. Nikulin
\paper The Igusa modular forms and ``the simplest''
Lorentzian Kac--Moody algebras
\jour Matem. Sbornik 
\yr 1996 \vol 187 
\pages 1601--1643 
\endref

\ref
\key GN3
\by V.A. Gritsenko, V.V. Nikulin
\paper Automorphic forms and Lorentzian Kac-Moody algebras.
Part I 
\jour International  J. of Mathem.
\vol 9
\yr 1998
\pages 153--199
\endref

\ref
\key GN4 
\by V.A. Gritsenko, V.V. Nikulin
\paper Automorphic forms and Lorentzian Kac-Moody algebras.
Part II 
\jour International  J. of Mathem.
\vol 9
\yr 1998
\pages 201--275
\endref


\ref
\key H1
\by F. Hirzebruch
\paper Elliptic genera of level N for complex manifolds
\inbook Differential geometrical Methods
in Theoretical Physics 
\eds K. Bleuler, M. Werner
\publ Kluwer Acad. Publ.
\yr 1988
\pages 37--63
\moreref Appendix III to [HBJ]
\endref

\ref
\key H2
\by F. Hirzebruch
\paper Letter to V. Gritsenko from  11 August 1997
\endref


\ref
\key HBJ
\by F. Hirzebruch, T. Berger, R. Jung
\book Manifolds and Modular forms
\publ Aspects of Math. {\bf E20}
Vieweg--Verlag
\yr 1992
\endref

\ref
\key H\"o
\by G. H\"ohn
\paper  Komplex elliptische Geschlechter und 
$S^1$-\"aquivariante Kobordismustheorie
\jour
 Diplomarbeit
\yr 1991
\publaddr Bonn
\endref


\ref
\key K
\by I. Krichever
\paper Generalized elliptic genera and Baker--Akhiezer functions
\jour Mat. Zametki
\vol 47
\yr 1990
\pages 34--45
\endref

\ref
\key KYY
\by
T. Kawai, Y. Yamada, S.-K. Yang
\paper
Elliptic Genera
and N=2 Superconformal Field Theory
\jour
Nucl. Phys. 
\vol B414
\yr 1994
\pages 191-212
\endref

\ref
\key KM
\by
T. Kawai, K. Mohri
\paper
Geometry of $(0,2)$ Landau--Ginzburg orbifolds
\jour
Nucl. Phys. 
\vol B425
\yr 1994
\pages 191-216
\endref


\ref
\key L
\by P.S. Landweber Ed.
\book Elliptic Curves and Modular Forms in 
Algebraic Topology
\publ Springer-Verlag
\yr 1988
\endref 

\ref 
\key Li
\by K. Liu
\paper On elliptic genera and theta-functions
\jour Topology
\vol 35
\pages 617--640
\yr 1996
\endref


\ref
\key S
\by S. M. Salamon
\paper On the cohomology of K\"ahler and hyper-K\"ahler manifolds
\jour Topology
\vol 35 \yr 1996
\pages 137--155
\endref


\ref
\key SVW
\by
S. Sethi, C. Vafa, E. Witten
\paper Constraints on low-dimensional string
compactifications
\jour   Nucl. Phys. 
\vol 
B480
\yr 1996
\pages 213--224
\endref

\ref
\key Sk
\by N-P. Skoruppa
\paper Modular forms
\moreref Appendix I in [HBJ], pp. 121--161
\endref 


\ref 
\key T
\by J. T. Tate
\paper The arithmetic of elliptic curves
\jour Invent. Math.
\vol 23
\pages 179--206
\yr 1974
\endref

\ref
\key W1
\by E. Witten
\paper The index of the Dirac operator in loop space
\inbook Elliptic Curves and Modular Forms in 
Algebraic Topology
\publ Springer-Verlag
\yr 1988
\pages 161--181
\endref


\ref
\key W2
\by
E. Witten
\paper On the Landau--Ginzburg description of 
$N=2$ minimal models
\jour Int. J. Mod. Phys.
\vol A9
\yr 1994
\pages 4783
\endref



\endRefs
\enddocument 
\end